\documentclass[aps,twocolumn,nofootinbib,longbibliography]{revtex4-1}
\usepackage[babel]{csquotes}
\usepackage{graphicx}
\usepackage{amsmath,amssymb}
\usepackage[colorlinks,citecolor=blue,linkcolor=blue,urlcolor=blue]{hyperref}
\usepackage{mathrsfs}
\usepackage{enumerate}
\def\be{\begin{equation}}
\def\M{$\cal M$ }
\def\ee{\end{equation}}

\begin{document}

\title{Michelangelo's Stone: \\ an Argument against Platonism in Mathematics}

\author{Carlo Rovelli\vspace{.5em} }
\affiliation{\small
\mbox{CPT, Aix-Marseille Universit\'e, Universit\'e de Toulon, CNRS,} \\
 \mbox{Samy Maroun Research Center for Time, Space and the Quantum,}  \\ 
Case 907, F-13288 Marseille, France. \vspace{.3em}}
\date{\small\today}

\begin{abstract}
\noindent  
If there is a `platonic world' \M of mathematical facts, what does  \M  contain precisely? I observe that if  \M is too large, it is uninteresting, because the value is in the selection, not in the totality; if it is smaller and interesting, it is not independent of us. Both alternatives challenge mathematical platonism.  I suggest that the  universality of \emph{our} mathematics may be a prejudice hiding its contingency, and illustrate contingent aspects of classical geometry, arithmetic and linear algebra. 
\end{abstract}

\maketitle

Mathematical platonism \cite{Plato} is the view that mathematical reality exists by itself, independently from our own intellectual activities.\footnote{For an overview, see for instance \cite{Linnebo2013}.}  Many top level mathematicians hold this view dear, and express the sentiment that they do not ``construct" new mathematics, but  rather ``discover"   structures that already exist: real entities in a platonic mathematical world.\footnote{Contemporary mathematicians that have articulated this view in writing include Roger Penrose \cite{Penrose2005} and Alain Connes \cite{Connesa}.} Platonism is alternative to other views on the foundations of mathematics, such as reductionism, formalism, intuitionism, or the Aristotelian idea that mathematical entities exist, but they are embodied in the material world \cite{Gillies}. 

Here, I present a simple argument \emph{against} platonism in mathematics, which I have not found in the literature.  

The argument is based on posing a question.  Let us assume that a platonic world of mathematical entities and mathematical truths does indeed exist, and is independent from us. Let us call this world $\cal M$, for Mathematics.  The question I consider is: what is it reasonable to expect  \M to contain?  

I argue that even a superficial investigation of this question reduces the idea of the existence of $\cal M$ to something trivial, or contradictory.  

In particular, I argue that the attempt to restrict \M to ``natural and universal" structures is illusory: it is the perspectival error of mistaking ourselves as universal  (``English is clearly the most natural of languages".)  In particular, I point out contingent aspects of the two traditional areas of classical mathematics: geometry and arithmetic, and of a key tool of modern science: linear algebra. 

\subsubsection*{Michelangelo's stone and Borges' library}

Say we take a Platonic stance about math: in some appropriate sense, the mathematical world \M exists. The expressions ``to exist", ``to be real" and similar can have a variety of meanings and usages, and this is a big part of the issue, if not the main one.  But for the sake of the present argument I do not need to define them---nor, for that matter, platonism---precisely. The argument remains valid, I think, under most reasonable definitions and usages. 

So, what does \M include?  

Certainly \M includes all the beautiful mathematical theories that mathematicians have discovered so far. This is  its primary purpose.  It includes Pythagoras' theorem, the  classification of the Lie groups, the properties of prime numbers and so on.  It includes the real numbers, and Cantor's proof that they are ``more" than the integers, in the sense defined by Cantor, and both possible extensions of arithmetic: the one where there are infinities larger than the integers and smaller than the reals, and the one where there aren't. It contains game theory, and topos theory. It contains lots of stuff. 

But \M cannot include \emph{only} what mathematicians have discovered so far, because the point of platonism is precisely that what they \emph{will} discover tomorrow already exists in the platonic world \emph{today}.  It would not be kind towards future mathematicians to expect that \M contains just a bit more of what we have already done.   Obviously whoever takes a Platonic view must assume that in the platonic world of math there is \emph{much} more than what has been already discovered.  How much?

Certainly \M contains, say, all true (or at least all demonstrable) theorems about integer numbers. All possible true theorems about Euclidean geometry, including all those we have not yet discovered.  But there should be more than that, of course, as there is far more math than that. 

We can get a grasp on the content of \M from the axiomatic formulation of mathematics: given various  consistent sets ${\cal A}_1,{\cal A}_2,...$ of axioms, \M will include all true theorems following from each ${\cal A}_n$, all waiting to be discovered by us.  We can list many sets of interesting axioms, and imagine  the platonic world \M to be the ensemble of theorems these imply, all nicely ordered in families, according to the corresponding set of axioms.   

But this is still insufficient, because a good mathematician tomorrow could come out with a \emph{new} set of axioms and find new great mathematics, like the people who discovered non-commutative geometry, or those who defined $C^*$ algebras did. 

We are getting close to what the platonic world must contain.  Let us assume that a sufficiently universal language exist, say based on logic.  Then the platonic world \M is the ensemble of all theorems that follow from all (non contradictory) choices of axioms.    This is a good picture of what \M  could be.  We have found the content of the platonic world of math.

But something starts to be disturbing. The resulting \M is big, extremely big: it contains too much junk.  The large majority of coherent sets of axioms are totally irrelevant.   

Before discussing the problem with precision, a couple of similes can help to understand what is going on. \\

\emph{(i)} During the Italian Renaissance, Michelangelo Buonarroti, one of the greatest artists of all times,  said that a good sculptor does not create a statue: he simply ``takes it out" from the block of stone where the statue already lay hidden. A statue is already there, in its block of stone.  The artists must simply expose it, carving away the redundant stone \cite{Neret2004}.  The artist does not ``create" the statue: he ``finds" it. 

In a sense, this is true: a statue is just a subset of the grains of stone forming the original block. It suffices to take away the other grains, and the statue is taken out.  But the hard part of the game is of course to find out \emph{which} subset of grains of stone to leave there, and this, unfortunately, is not written on the stone.  It is selection that matters.  A block of stone already contained Michelangelo's \emph{Moses}, but it also contained virtually anything else --that is, all possible forms. The art of sculpture is to be able to determine, which, among this virtual infinity of forms, will talk to the rest of us as art does. 

Michelangelo's statement is evocative, maybe powerful, and perhaps speaks to his psychology and his greatness.  But it does not say much about the art of the sculptor.  The fact of including all possible statues does not confer to the stone the immense artistic value of all the possible beautiful statues it might contain, because the point of the art is the choice, not the collection. By itself, the stone is dull.\\

\emph{(ii)} The same story can be told about books.  Borges's famous library contained all possible books: all possible combinations of alphabet letters \cite{Borges2000}. Assuming a book is shorter than, say, a million letters, there are then more or less $30^{10^6}$ possible books, which is not even such a big number for a mathematician.  So, a writer does not really create a book: she simply ``finds" it, in the platonic library of all books. A particularly nice combinations of letters makes up, say, \emph{Moby-Dick}.  \emph{Moby-Dick} already existed in the platonic space of books: Melville didn't create  \emph{Moby-Dick}, he just discovered  \emph{Moby-Dick}...    

Like Michelangelo's stone, Borges's library is void of any interest: it has no content, because the value is in the choice, not in the totality of the alternatives. \\

The Platonic world of mathematics \M defined above is similar to Michelangelo's block of stone, or Borges library, or Hegel's ``night in which all cows are black"\footnote{Hegel utilized this Yiddish saying to ridicule Schelling's notion of Absolute, meaning that --like mathematical platonism-- this included too much and was too undifferentiated, to be of any relevance \cite{Hegel}.}: a mere featureless vastness, without value because the value is in the choice, not in the totality of the possibilities.   Similarly, science can be said to be nothing else than the denotation of a subset of David Lewis's possible worlds \cite{Lewis2001}: those respecting certain laws we have found. But Lewis's totality of all possible world is not science, because the value of science  is in the restriction, not in the totality.  

Mathematics may be called an ``ensemble of tautologies", in the sense of Wittgenstein. But it is not the ensemble of \emph{all} tautologies: these are too many and their ensemble is uninteresting and irrelevant. Mathematics is about recognizing the ``interesting" ones.  Mathematics may be the investigation of structures. But it is not the list of \emph{all} possible structures: these are too many and their ensemble is uninteresting.  If the world of mathematics was identified with the platonic world \M defined above, we could program a computer to slowly unravel it entirely, by listing all possible axioms and systematically applying all possible transformation rules to derive all possible theorems.  But we do not even think of doing so. Why? Because what we call mathematics is an infinitesimal subset of the huge world $\cal M$ defined above: it is the tiny subset which is of interest for us.  Mathematics is about studying the ``interesting" structures.  

So, the problem becomes: what does ``interesting" mean?

\subsubsection*{Interest is in the eye of the interested}

Can we restrict the definition of \M to the \emph{interesting} subset? Of course we can, but interest is in the eyes of a subject.   A statue is a subset of the stone which is worthwhile, \emph{for us}. A particular combination of letters is a good book, \emph{for us}.    What is it that makes certain set of axioms defining certain mathematical objects, and certain theorems, interesting? 

There are different possible answer to this question, but they all make explicit or implicit reference to features of ourselves, our mind, our contingent environment, or the physical structure our world happens to have.  

This fact is pretty evident as far as art or literature are concerned. Does it hold for mathematics as well?  Hasn't mathematics precisely that universality feel that is at the root of platonism? 

Shouldn't we expect ---as often claimed--- \emph{any} other intelligent entity of the cosmos to come out with the \emph{same} ``interesting" mathematics as us? 

The question is crucial for mathematical platonism, because platonism is the thesis that mathematical entities and truths form a world which exists \emph{independently from us}. If what we call mathematics ends up depending heavily on ourselves or the specific features of our world, platonism looses its meaning.

I present below some considerations that indicate that the claimed universality of mathematics is a parochial prejudice.  These are based on the concrete examples provided by the chapters of mathematics that have most commonly been indicated as universal. 

\subsubsection*{The geometry of a sphere}

Euclidean geometry has been among the  first pieces of mathematics to be formalized.  Euclid's celebrated text, the ``Elements" \cite{Euclid1997}, where Euclidean geometry is beautifully developed, has been the ideal reference for all mathematical texts. Euclidean geometry describes space. It captures our intuition about the structure of space. It has  applications to virtually all fields of science, technology and engineering. Pythagoras' theorem, which is at its core, is a familiar icon.  

It is difficult to imagine something having a more ``universal" flavor than euclidean geometry. What could be contingent or accidental about it? What part do we humans have in singling it out?   Wouldn't any intelligent entity developing anywhere in the universe come out with this same mathematics?

I maintain the answer is negative.  To see why, let me start by recalling that, as is well known, Euclidean geometry was developed by Greek mathematicians mostly living in Egypt during the Hellenistic period, building on Egyptians techniques for measuring the land. These were important because of the Nile's floods cancelling borders between private land parcels. The very name of the game, ``geometry", means  ``measurement of the land" in Greek.  Two-dimensional Euclidean geometry describes, in particular, the mathematical structure formed by the land. 

But: does it? 

Well, the Earth is more a sphere than a plane. Its surface is better described by the geometry of a sphere, than by two-dimensional (2d) Euclidean geometry.  It is an accidental fact that Egypt happens to be small compared to the size of the Earth. The radius of the Earth is around 6,000 Kilometers.  The size of Egypt is of the order of 1,000 Kilometers. Thus, the scale of the Earth is more than 6 times larger than the scale of Egypt. Disregarding the sphericity of the Earth is an approximation, which is viable when dealing with the geometry of Egypt and becomes better and better as the region considered is smaller. As a practical matter, 2d Euclidean geometry is useful, but it is a decent approximation that works \emph{only} because of the smallness of the size of Egypt.  Intelligent beings living on a planet just a bit smaller than ours \cite{DeSaint-Exupery1945}, would have easily detected the effects of the curvature of the planet's surface. They would not have developed 2d Euclidean geometry.  

One may object that this is true for 2d, but not for 3d geometry. The geometry of the surface of a sphere can after all be obtained from Euclidean 3d geometry.  But the objection has no teeth: we have learned with general relativity that spacetime is curved and Euclidean geometry is just an approximation also as far as 3d physical space is concerned. Intelligent beings living on a region of stronger spacetime curvature would have no reason to start mathematics from Euclidean geometry.\footnote{It is well known that Kant was mistaken in his deduction that the Euclidean geometry of physical space is true a priori \cite{Kant}. But even Wittgenstein bordered on mistake in dangerously appearing to assume a unique possible set of  laws of geometry for anything spatial: ``We could present spatially an atomic fact which contradicted the laws of physics, but not one which contradicted the laws of geometry".  Tractatus, Proposition 3.0321 \cite{Wittgenstein}.}

\begin{figure}
\center
\raisebox{5mm}{\includegraphics[width=5cm]{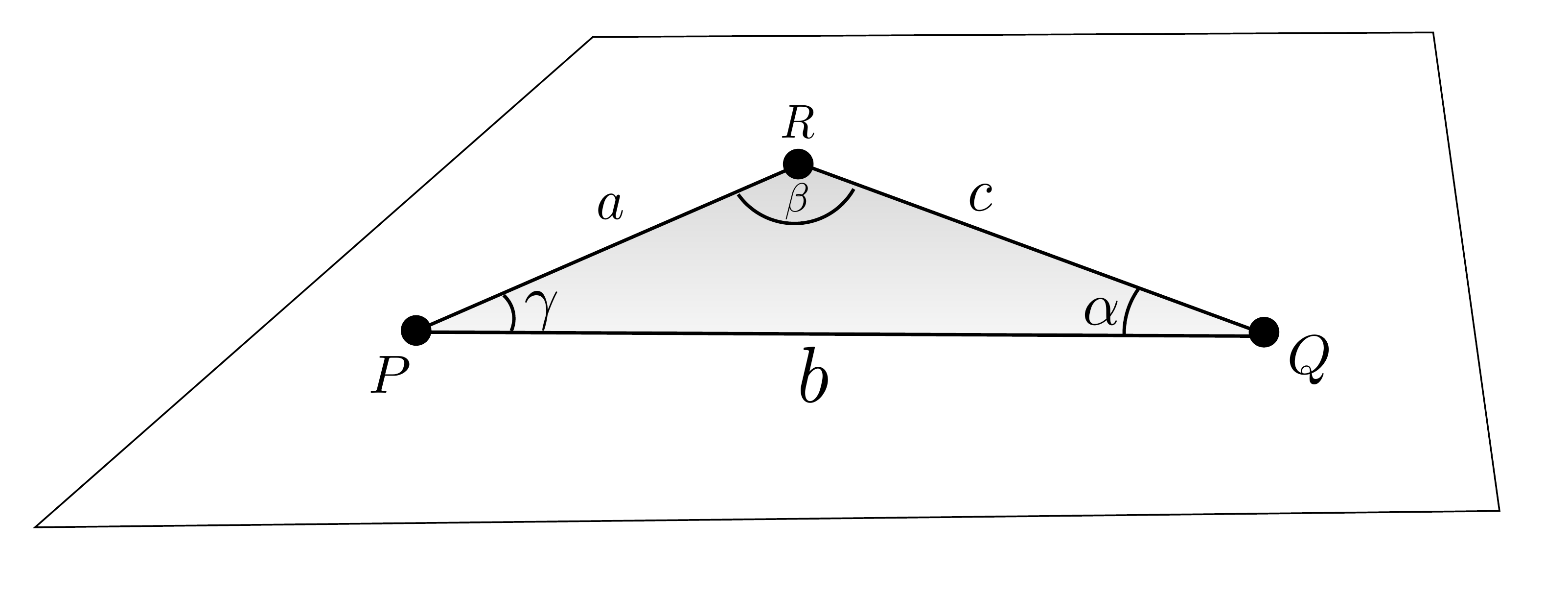}}\ \ \ 
\includegraphics[width=3cm]{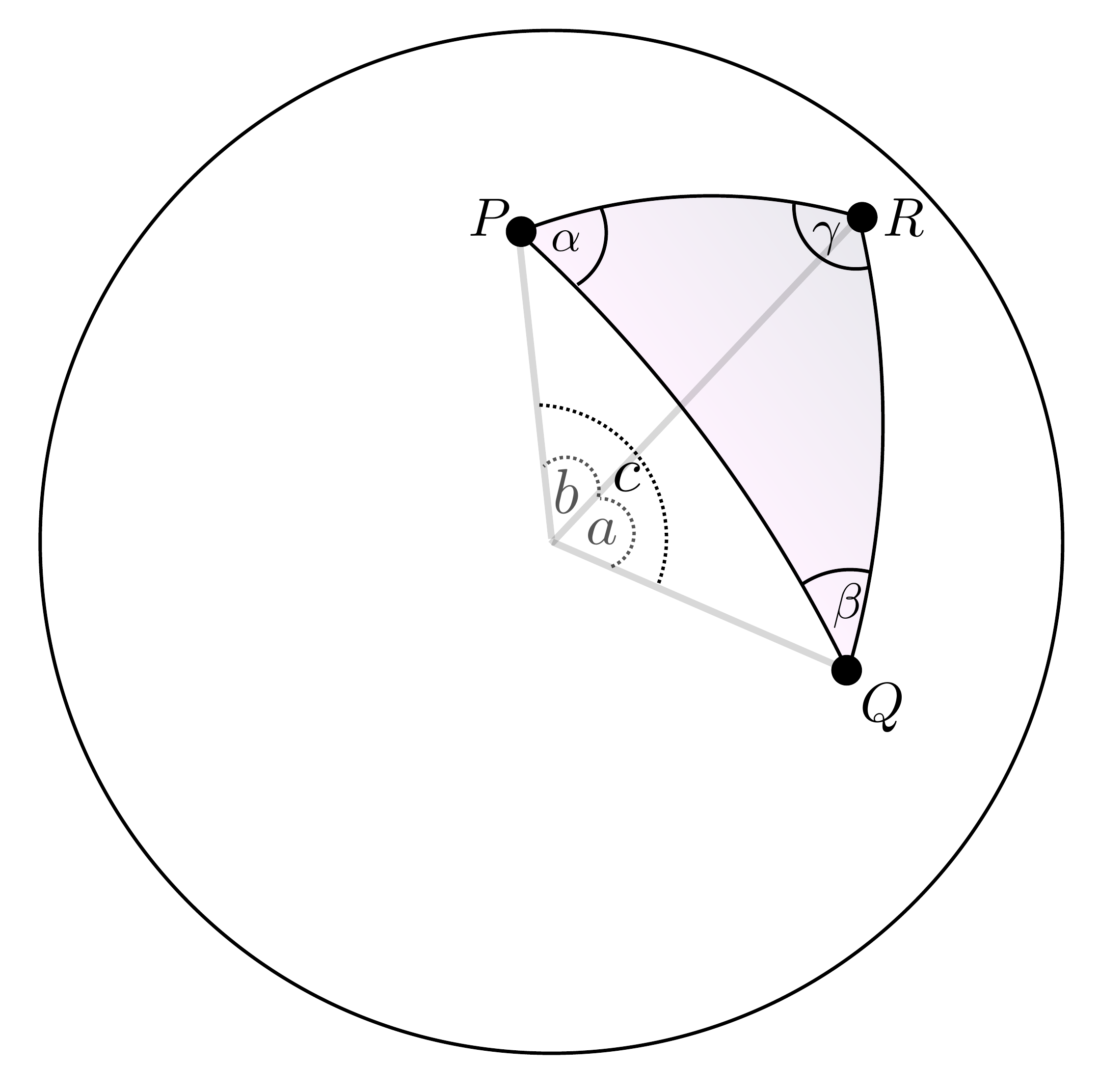}
\caption{Flat and spherical triangles.}
\end{figure}

A more substantial objection is that 2d euclidean geometry is \emph{simpler} and more ``natural" than curved geometry.  It is intuitively grasped by our intellect, and mathematics describes this  intuition about space. Its simplicity and intuitive aspect are the reasons for its universal nature. Euclidean geometry is therefore universal in this sense. I show below that  this objection is equally ill founded: the geometry of a sphere is definitely simpler and more elegant than the geometry of a plane. 

Indeed, there is a branch of mathematics, called 2d ``spherical" geometry, which describes directly the (intrinsic) geometry of a sphere.  This is the mathematics that the Greeks would have developed had the Earth been sufficiently small to detect the effects of the Earth's surface curvature on the Egyptians fields.  Perhaps quite surprisingly for many, spherical geometry is far simpler and ``more elegant" than  Euclidean geometry.  I illustrate this statement with a few examples below, without, of course, going into a full exposition of spherical geometry (see for instance \cite{Todhunter1886,Harris1998}). 

Consider the theory of triangles: the familiar part of geometry we study early at school.  In Euclidean geometry, a triangle has three sides, with lengths, $a$, $b$ and $c$, and three angles $\alpha$, $\beta$ and $\gamma$ (Figure 1). We measure angles with pure numbers, so, $\alpha$, $\beta$ and $\gamma$ are taken to be numbers with value between $0$ and $\pi$.  Measuring with numbers the length of the sides is a more complicated business.  Since there is no standard of length in Euclidean geometry, we must either resort to talk only about ratios between lengths (as the ancients preferred), or to    choose arbitrarily a segment once and for all, use it as our ``unit of measure", and characterize the length of each side of the triangle by the number which gives its ratio to the unit (as the moderns prefer).

\begin{figure}
\includegraphics[width=4cm]{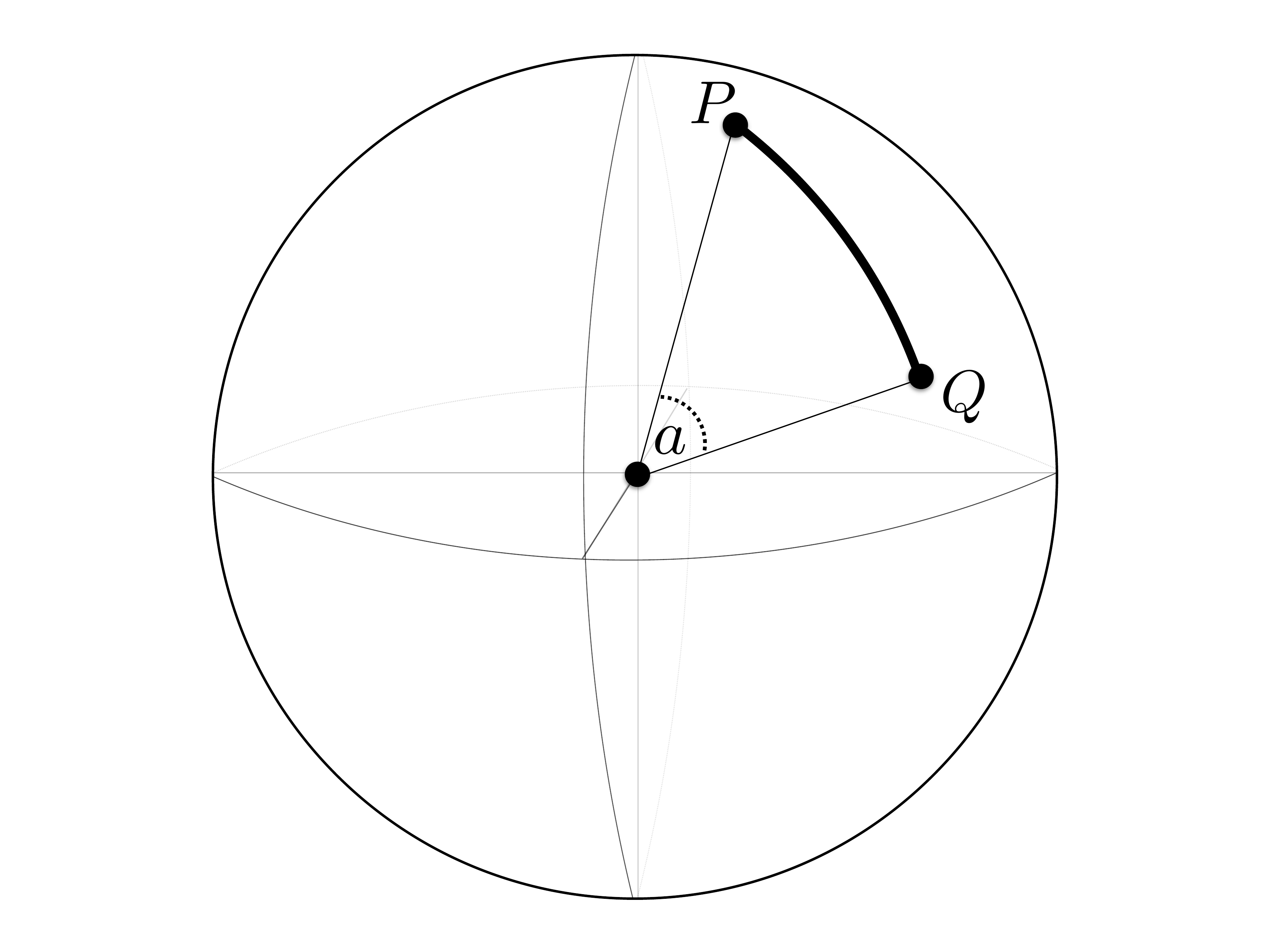}
\caption{Two points on a sphere determine an arc, whose size is measured by the angle it subtends, or equivalently, intrinsically, by its ratio to an equator.}
\end{figure}

All this simplifies dramatically in spherical geometry: here there is a preferred scale, the length of the equator. The length of an arc (the shortest path between two points) is naturally measured by the ratio to it. Equivalently (if we immerge the sphere in space) by the angle subtended to the arc. Therefore the length of the side of a triangle ($a$, $b$, $c$) is an angle \emph{as well}. See Figure 2. Compare then the theories of triangles in the two geometries (Figure 1):

{\em Euclidean geometry:}
\vspace{-2mm}
\begin{enumerate}[(i)]\addtolength{\itemsep}{-2mm}
\item Two triangles are equal if
they have equal sides, or if one side and two angles are equal. 
\item The area of the triangle is \\ $A=\frac14\sqrt{2a^2b^2+2a^2c^2+2b^2c^2-a^4-b^4-c^4}$.
\item For right triangles: $a^2+b^2=c^2$.
\end{enumerate}

{\em Spherical geometry: }
\vspace{-2mm}
\begin{enumerate}[(i)]\addtolength{\itemsep}{-2mm}
\item Triangles with  same sides, or same angles, are equal.
\item The area of a triangle is  $A=\alpha+\beta+\gamma-\pi$. 
\item For right triangles: $\cos c=\cos a \cos b$.
\end{enumerate}

Even a cursory look at these results reveals the greater simplicity of spherical geometry. Indeed, spherical geometry has a far more ``universal" flavor than Euclidean geometry.  

Euclidean geometry can be be obtained from it as a limiting case: it is the geometry of figures that are much smaller than the curvature radius. In this case $a$, $b$ and $c$ are all much smaller than $\pi$. Their cosine is well approximated by $\cos\theta\sim 1-\frac12\theta^2$ and the last formula reduces to  Pythagoras' theorem in the first approximation.  Far from being a structural property of the very nature of space, Pythagoras' theorem is only a first order approximation, valid in a limit case of a much simpler and cleaner mathematics: 2d spherical geometry. 

There are many other beautiful and natural results in spherical geometry, which I shall not report here. They extend to the 3d case: the intrinsic geometry of a 3-sphere.  A 3-sphere is a far more reasonable structure than the infinite Euclidean space: it is the finite homogeneous three-dimensional metric space without boundaries.  The geometry may well be the large scale geometry of our universe \cite{Einstein:1917ce}.\footnote{Cosmological measurements indicate that \emph{spacetime} is curved, but have so far failed to detected a large scale cosmological curvature of \emph{space}. This of course does \emph{not} imply that the universe is flat \cite{Ellis:2005pv}, for the same reason for which the failure to detect curvature on the fields of Egypt did not imply that that the Earth was flat. It only shows that the universe is big.  Current measurements indicate that the radius of the Universe should be at least ten time larger than the portion of the Universe we see \cite{Collaboration2009}.  A ratio, by the way, quite similar to the Egyptian case.}  It shape is counterintuitive for many of us, schooled in Euclid. But it was not so for Dante Alighieri, who  did not study Euclid at school: the topology of the universe he describes in his poem is precisely that of a 3-sphere \cite{Peterson:1979fk}. See Figure 3. What is ``intuitive" changes with history. 

\begin{figure}[t]
\begin{center}
\includegraphics[width=2in]{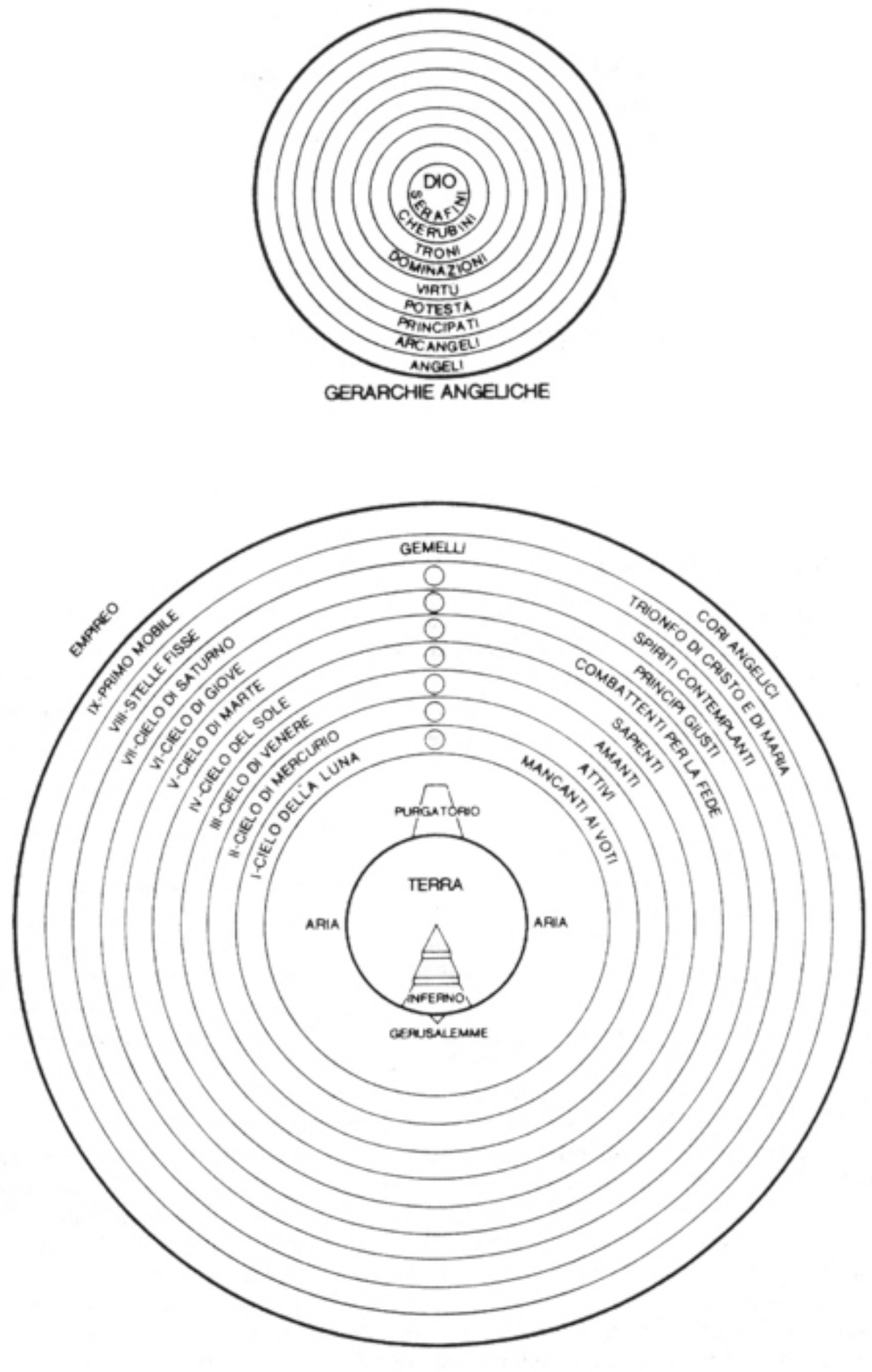}\vspace{-5mm}
\caption{Dante's universe: the Aristotelian spherical universe is surrounded by another similar spherical space, inhabited by God and Angel's spheres.  The two spheres together form a three-sphere.}
\label{Dante}
\end{center}
\end{figure}

These considerations indicate that the reason Euclidean geometry has played such a major role in the foundation of mathematics is not because of its universality and independence from our contingent situation. It is the opposite: Euclidean geometry is of value to us just because it describes---not even very well---the accidental properties of the region we happen to inhabit.  Inhabitants of a different region of the universe---a smaller planet, or a region with high space curvature---would likely fail to consider euclidean geometry interesting mathematics. For them, Euclidean geometry could be a uninteresting and cumbersome limiting case.

\subsubsection*{Linear algebra}

Every physicist, mathematician or engineer learns linear algebra and uses it heavily.  Linear algebra, namely the algebra of vectors, matrices, linear transformations and so on, is the algebra of linear spaces, and since virtually everything is linear in a first approximation, linear algebra is ubiquitous. It is difficult to resist its simplicity, beauty and generality when studying it, usually in the early years at college. 

Furthermore, today, we find linear algebra at the very foundations of physics, because it is the language of quantum theory.  In the landmark paper that originated quantum theory \cite{Heisenberg1925a}, Werner Heisenberg understood that physical quantities are better described by \emph{matrices} and used the multiplication of matrices (core of linear algebra) to successfully compute the properties of quantum particles.  Shortly later, Paul Dirac wrote his masterpiece book \cite{Dirac1930}, where quantum mechanics is entirely expressed in terms of linear algebra (linear operators, eigenvalues, eigenvectors...).

It would therefore seem natural to formulate the hypothesis that any minimally advanced civilization would  discover linear algebra very early and start using it heavily.

But it is not the case.  In fact, we have a remarkable counterexample: a civilization that has developed for millennia without developing linear algebra---ours. 

When Heisenberg wrote his famous paper he did not know linear algebra. He had no idea of what a matrix is, and had never previously learned the algorithm for multiplying matrices. He made it up in his effort to understand a puzzling aspect of the physical world.  This is pretty evident from his paper.  Dirac, in his book, is basically inventing linear algebra in the highly non-rigorous manner of a physicist. After having constructed it and tested its power to describe \emph{our} world, linear algebra appears natural to us.  But it didn't appear so for generations of previous mathematicians.  

Which tiny piece of \M turns out to be interesting for us,  which parts turns out to be ``mathematics" is far from obvious and universal. It is largely contingent. 

\subsubsection*{Arithmetic and identity}

The last example I discuss is given by the natural numbers (1, 2, 3, ...), which form the basis of arithmetic, the other half of classical mathematics.   Natural numbers seem very natural indeed. There is evidence that the brain is pre-wired to be able to count, and do elementary arithmetic with small numbers \cite{Vallortigara}.  Why so?  Because our world appears to be naturally organized in terms of things that can be counted.   But is this a feature of reality at large, of any possible world, or is it just a special feature of this little corner of the universe we inhabit and perceive?   

I suspect the second answer to be the right one.  The notion of individual ``object" is notoriously slippery, and objects need to have rare and peculiar properties in order to be countable. How many clouds are there in the sky? How many mountains in the Alps? How many coves along the coast of England?  How many waves in a lake? How many clods in my garden?  These are all very ill-defined questions. 

To make the point, imagine  some form of intelligence evolved on Jupiter, or a planet similar to Jupiter. Jupiter is fluid, not solid.  This does not prevent it from developing complex structures: fluids develop complex structures, as their chemical composition, state of motion, and so on, can change \emph{continuously} from point to point and from time to time, and their dynamics is governed by rich nonlinear equations. Furthermore, they interact with  magnetic  and electric fields, which vary continuously in space and time as well.  Imagine that in such a huge (Jupiter is much larger than Earth's) Jovian environment, complex structures develop to the point to be conscious and to be able to do some math.  After all, it has happened on Earth, so it shouldn't be so improbable for something like this to happen on an entirely fluid planet as well.  Would this math include counting, that is, arithmetic? 

Why should it?  There is nothing to count in a completely fluid environment.  (Let's also say that our Jovian planet's atmosphere is so thick that one cannot see and count the stars, and that the rotation and revolution periods are equal, as for our Moon, and there are neither days nor years.)  The math needed by this fluid intelligence would presumably include some sort of geometry, real numbers, field theory, differential equations..., all this  could develop using only geometry, without ever considering this funny operation which is enumerating individual things one by one.   

The notion of ``one thing", or ``one object", the notions themselves of unit and identity, are useful for us living in an environment where there happen to be stones, gazelles, trees, and friends that can be counted. The fluid intelligence diffused over the Jupiter-like planet, could have developed mathematics without ever thinking about natural numbers. These would not be of interest for her. 

I may risk being more speculative here. The development of the ability to count may be connected to the fact that life evolved on Earth  in a peculiar form characterized by the existence of  ``individuals".   There is no reason an intelligence capable to do math should take this form. In fact, the reason counting appears so natural to us may be that we are a species formed by interacting individuals, each realizing a notion of identity, or unit.   What is clearly made by  units is a group of interacting primates, not the world.  The archetypical identities are my friends in the group.\footnote{This might be why ancient humans attributed human-like mental life to animals, trees and stones: they were perhaps utilizing mental circuits developed to deal with one another --within the primate group-- extending them to deal also with animals, trees and stones.}   

Modern physics is intriguingly ambiguous about countable entities. On the one hand, a major discovery of the XX century has been that at the elementary level nature is entirely described by field theory.  Fields  vary continuously in space and time.  There is little to count, in the field picture of the world.  On the other hand, quantum mechanics has injected a robust dose of discreteness in fundamental physics: because of quantum theory, fields have particle-like properties and particles are quintessentially countable objects. In any introductory quantum field theory course,  students meet an important operator, the number operator $N$, whose eigenvalues are the natural numbers and whose physical interpretation is counting particles \cite{Itzykson2006}.  Perhaps our fluid Jovian intelligence would finally get to develop arithmetic when figuring out quantum field theory...   

But notice that what moderns physics says about what is countable in the world has no bearing on the universality of mathematics: at most, it points out which parts of \M are interesting because they happen to describe \emph{this} world.  

\subsubsection*{Conclusion}

In the light of these consideration, let us look back at the development of our own mathematics. Why has mathematics developed at first, and for such a long time, along two parallel lines: geometry and arithmetic?  The answer begins to clarify: because these two branches of mathematics are of  value for creatures like us, who instinctively count friends, enemies and sheep, and who need to measure, approximately, a nearly flat earth in a nearly flat region of physical space.   In other words, this mathematics is of interest to us because it reflects very contingent interests of ours.    Out of the immense vastness of $\cal M$, the dull platonic space of all possible structures, we have carved out, like Michelangelo, a couple of shapes that talk to us.  From the immense vastness of M, the dull platonic space of all possible structures, we have carved out, like Michelangelo, a couple of shapes that speak to us.

There is no reason to assume that the mathematics that has developed later escapes this contingency.
To the contrary, the continuous re-foundations and the constant re-organization of the global structure of mathematics testify to its non-systematic and non-universal global structure. Geometry, arithmetic, algebra, analysis, set theory, logic, category theory, and --recently-- topos theory \cite{Caramello2010} have all been considered for playing a foundational role in mathematics.  Far from being stable and  universal, our mathematics is a fluttering butterfly, which follows the fancies of inconstant creatures.  Its theorems are solid, of course; but selecting what represents an interesting theorem is a highly subjective matter. 

It is the carving out, the selection, out of a dull and undifferentiated $\cal M$, of a subset which is useful \emph{to us}, interesting \emph{for us}, beautiful and simple in \emph{our} eyes, it is, in other words, something strictly related to what \emph{we} are, that makes up what we call mathematics.   

The idea that the mathematics that we find valuable forms a Platonic world \emph{fully independent from us} is like the idea of an Entity that created the heavens and the earth, and happens to very much resemble my grandfather. 

\centerline{---}

Thanks to Hal Haggard and Andrea Tchertkoff for a careful reading of the manuscript and comments. 

\providecommand{\href}[2]{#2}\begingroup\raggedright\endgroup

\end{document}